\documentclass[12pt,reqno]{amsart}
\usepackage{amsmath}
\usepackage{amsthm}
\usepackage{amsfonts}
\usepackage{mathrsfs}
\usepackage{paralist}
\usepackage{cite}
\usepackage[hypertex]{hyperref}
\newtheorem{theorem}{Theorem}[section]
\newtheorem{prop}[theorem]{Proposition}
\newtheorem{problem}[theorem]{Problem}
\theoremstyle{definition}
\newtheorem{definition}[theorem]{Definition}
\theoremstyle{remark}
\newtheorem{example}[theorem]{Example}
\newtheorem{remark}[theorem]{Remark}
\newtheorem*{ackn}{Acknowledgements}
\DeclareMathOperator{\Hom}{Hom}
\DeclareMathOperator{\Pol}{Pol}
\newcommand*{\Pfree}{\mathop{\widehat *}}
\newcommand*{\cHol}{\mathop{\mathit{Hol}}}
\newcommand*{\op}{\mathrm{op}}
\newcommand*{\reg}{\mathrm{reg}}

\newcommand*{\wh}{\widehat}
\newcommand*{\la}{\langle}
\newcommand*{\ra}{\rangle}
\newcommand{\Alg}{\mathsf{Alg}}
\newcommand{\ArM}{\mathsf{AM}}
\newcommand*{\CC}{\mathbb C}
\newcommand*{\N}{\mathbb N}
\newcommand*{\Z}{\mathbb Z}
\newcommand*{\R}{\mathbb R}
\newcommand*{\DD}{\mathbb D}
\newcommand*{\BB}{\mathbb B}
\newcommand*{\cO}{\mathscr O}
\newcommand*{\cF}{\mathscr F}
\newcommand*{\fF}{\mathfrak F}
\newcommand*{\ccH}{\mathcal H}
\newcommand*{\cN}{\mathscr N}
\newcommand*{\cB}{\mathscr B}
\newcommand*{\cU}{\mathscr U}

\newcommand*{\defo}{\mathrm{def}}
\newcommand*{\eps}{\varepsilon}
\newcommand*{\ol}{\overline}

\begin{document}
\title[Quantum polydisk and quantum ball]{Quantum polydisk, quantum ball, and a $q$-analog\\
of Poincar\'e's theorem}
\author{A. Yu. Pirkovskii}
\address{Faculty of Mathematics\\
National Research University ``Higher School of Economics''\\
7~Vavilova, 117312 Moscow, Russia}
\email{aupirkovskii@hse.ru, pirkosha@gmail.com}
\thanks{This work was supported by the RFBR grant no. 12-01-00577.}
\date{}

\begin{abstract}
The classical Poincar\'e theorem (1907) asserts that the polydisk $\DD^n$
and the ball $\BB^n$ in $\CC^n$ are not biholomorphically equivalent for $n\ge 2$.
Equivalently, this means that the Fr\'echet algebras $\cO(\DD^n)$ and
$\cO(\BB^n)$ of holomorphic functions are not topologically isomorphic. Our goal
is to prove a noncommutative version of the above result.
Given $q\in\CC\setminus\{ 0\}$, we define two noncommutative power series algebras
$\cO_q(\DD^n)$ and $\cO_q(\BB^n)$, which can be viewed as $q$-analogs of
$\cO(\DD^n)$ and $\cO(\BB^n)$, respectively.
Both $\cO_q(\DD^n)$ and $\cO_q(\BB^n)$ are the completions of the algebraic
quantum affine space $\cO_q^\reg(\CC^n)$ w.r.t. certain families of seminorms.
In the case where $0<q<1$, the algebra $\cO_q(\BB^n)$ admits an equivalent definition
related to L.~L.~Vaksman's algebra $C_q(\bar\BB^n)$ of continuous functions
on the closed quantum ball. We show that both $\cO_q(\DD^n)$ and $\cO_q(\BB^n)$
can be interpreted as Fr\'echet algebra deformations (in a suitable sense) of $\cO(\DD^n)$ and
$\cO(\BB^n)$, respectively. Our main result is that $\cO_q(\DD^n)$ and $\cO_q(\BB^n)$
are not isomorphic if $n\ge 2$ and $|q|=1$, but are isomorphic if $|q|\ne 1$.
\end{abstract}

\maketitle

\section{Introduction}
\label{sec:intro}
Noncommutative geometry is a vast and rapidly growing subject consisting
of a number of different branches (noncommutative algebraic geometry,
noncommutative differential geometry, noncommutative topology,
noncommutative measure theory, etc.). Each of these branches has its own
objects of study and its own methods. Nevertheless, all of them share
the common unifying ``philosophy'' that some classical constructions and
results known from various fields of geometry and topology can
be successfully applied to noncommutative objects, which, at the first
glance, have nothing to do with geometry.

The subject of the present paper can be characterized as ``noncommutative complex
analysis'', or ``noncommutative complex analytic geometry''. At the moment,
this theory is much less developed than any of the above-mentioned parts of
noncommutative geometry. However, a number of important results have been
obtained in this field during the last decade. First of all, let us mention the tremendous work
done by L.~L.~Vaksman's school
(see, e.g., \cite{Vaks_int_disk,Vaks_Dolb,Vaks_max,Vaks_sborn,Vaks_book}
and references therein), which eventually resulted in the creation of
the general theory of quantum bounded symmetric domains.
A more operator-theoretic aspect of the subject is reflected in
the papers by K.~R.~Davidson, D.~R.~Pitts,
E. G. Katsoulis, C. Ramsey, and O.~Shalit
\cite{Dav_FSG,Dav_Topl,Dav_Pop,Dav_dil,Dav_bihol,Dav_isoprobl},
F.~H.~Szafraniec \cite{Szfr_mult}, and G.~Popescu
\cite{Pop_disk,Pop_freeprod,Pop_holball,Pop_var_2,Pop_aut,Pop_ncdom,Pop_class,Pop_jntsim,Pop_freebihol1,Pop_freebihol2}.
An algebraic point of view is adopted by A.~Polishchuk and
A.~Schwarz \cite{Pol_Schw,Pol_clbnd,Pol_expmap,Pol_qsicoh},
M.~Khalkhali, G.~Landi, W.~D.~van Suijlekom, and A.~Moatadelro \cite{KLvS,KM1,KM2},
E.~Beggs and S.~P.~Smith \cite{Beggs_Smith}, R. \'O Buachalla \cite{Buach_qbnd,Buach_hmg}.
All this shows that noncommutative
complex analysis and noncommutative complex geometry will hopefully reach
their bloom period fairly soon.

Our primary goal is to prove a noncommutative
version of a classical theorem by Poincar\'e \cite{Poincare},
which asserts that the polydisk $\DD^n$
and the ball $\BB^n$ in $\CC^n$ are not biholomorphically equivalent for $n\ge 2$.
This result is often mentioned in textbooks as one of the first results in
function theory of several complex variables (see, e.g., \cite{Narasim,Range}).
To understand how a noncommutative version of Poincar\'e's theorem should look like,
let us recall O.~Forster's important result \cite{Forster} whose informal meaning is that
all essential information about a domain of holomorphy $D$ in $\CC^n$ is contained in the
algebra $\cO(D)$ of holomorphic functions on $D$. To be more precise, Forster's theorem
states that the functor
\begin{align*}
\Bigl\{ \text{ Stein spaces }\Bigr\} &\to \Bigl\{ \text{ Fr\'echet algebras }\Bigr\},\\
(X,\cO_X) &\mapsto \cO(X),
\end{align*}
is fully faithful.
As a consequence, two domains of holomorphy $D_1$ and $D_2$ in $\CC^n$ are
biholomorphically equivalent if and only if the algebras $\cO(D_1)$ and $\cO(D_2)$
are topologically isomorphic. Thus Poincar\'e's theorem is equivalent to the assertion that
the algebras $\cO(\DD^n)$ and $\cO(\BB^n)$ are not topologically isomorphic for $n\ge 2$.

Now we can explain our plan in more detail. Let $\CC^\times=\CC\setminus\{ 0\}$.
For each $q\in\CC^\times$
we define $q$-analogs of the algebras $\cO(\DD^n)$ and
$\cO(\BB^n)$ to be the completions of the algebraic
quantum affine space $\cO_q^\reg(\CC^n)$ with respect to certain families of seminorms.
Both the resulting algebras $\cO_q(\DD^n)$ and $\cO_q(\BB^n)$ were
introduced earlier in \cite{Pir_Shlmn} (see also \cite{Pir_HFG} for a more detailed
treatment of $\cO_q(\DD^n)$), but the definition of $\cO_q(\BB^n)$ was given only
for $0<q<1$. Here we propose a different approach to $\cO_q(\BB^n)$ which makes sense
for all $q\in\CC^\times$ and which is equivalent to the approach of \cite{Pir_Shlmn}
in the case where $0<q<1$. To justify our definitions, we show that the Fr\'echet
algebra families $\{ \cO_q(\DD^n) : q\in\CC^\times\}$
and $\{ \cO_q(\BB^n) : q\in\CC^\times\}$ can be arranged into Fr\'echet algebra bundles
over $\CC^\times$, generalizing thereby our earlier result from \cite{Pir_Shlmn}.
Our main result, i.e., a $q$-analog of Poincar\'e's theorem, states
that $\cO_q(\DD^n)$ and $\cO_q(\BB^n)$
are not topologically isomorphic if $n\ge 2$ and $|q|=1$. On the other hand,
we show that they are topologically isomorphic if $|q|\ne 1$.

This paper is mostly a survey. The proofs are either sketched or omitted. We plan to present
the details elsewhere.

\section{Preliminaries}
We shall work over the field $\CC$ of complex numbers. All algebras are assumed to
be associative and unital, and all algebra homomorphisms are assumed to be unital
(i.e., to preserve identity elements).
By a {\em Fr\'echet algebra} we mean a complete metrizable locally convex
algebra (i.e., a topological algebra whose underlying space
is a Fr\'echet space). A {\em locally $m$-convex algebra} \cite{Michael} is a topological
algebra $A$ whose topology can be defined by a family of submultiplicative
seminorms (i.e., seminorms $\|\cdot\|$ satisfying $\| ab\|\le \| a\| \| b\|$
for all $a,b\in A$). A complete locally $m$-convex algebra is called
an {\em Arens-Michael algebra} \cite{X2}. The algebra of holomorphic functions
on a complex manifold $X$ will be denoted by $\cO(X)$. Recall that $\cO(X)$ is
a Fr\'echet-Arens-Michael algebra with respect to the topology of uniform convergence
on compact subsets of $X$.

\section{Quantum affine space}
Let $q\in\CC^\times$.
Recall that the algebra $\cO_q^\reg(\CC^n)$ {\em of regular functions on the
quantum affine $n$-space} (see, e.g., \cite{Br_Good}) is generated by
$n$ elements $x_1,\ldots ,x_n$ subject to the relations $x_i x_j=qx_j x_i$ for all $i<j$.
If $q=1$, then $\cO_q^\reg(\CC^n)$ is nothing but the
polynomial algebra $\CC[x_1,\ldots ,x_n]=\cO^\reg(\CC^n)$.
Of course, $\cO_q^\reg(\CC^n)$ is noncommutative unless $q=1$, but the monomials
$x^k=x_1^{k_1}\cdots x_n^{k_n}\; (k\in\Z_+^n)$ still
form a basis of $\cO_q^\reg(\CC^n)$. Thus $\cO_q^\reg(\CC^n)$ may be viewed as
a ``deformed'' polynomial algebra.

The algebras $\cO_q(\DD^n_r)$ and $\cO_q(\BB^n_r)$ that we are going to define
will be certain completions of $\cO_q^\reg(\CC^n)$. There are many nonequivalent ways
to complete this algebra, but, among all the completions, there is a universal one.
Recall that the {\em Arens-Michael envelope}, $\wh{A}$, of an algebra $A$
is the completion of $A$ with respect to the family of all submultiplicative
seminorms on $A$.
The Arens-Michael envelope has the universal property that, for each Arens-Michael
algebra $B$, there is a $1$-$1$ correspondence
\[
\Hom_\Alg(A,B)\cong\Hom_\ArM(\wh{A},B),
\]
where $\Alg$ is the category of algebras and $\ArM$ is the category of
Arens-Michael algebras. Moreover, the assignment $A\mapsto\wh{A}$ is a functor
from $\Alg$ to $\ArM$, and this functor is left adjoint to the forgetful functor
$\ArM\to\Alg$.

Arens-Michael envelopes were introduced by J.~L.~Taylor \cite{T1} under the name of
``completed locally $m$-convex envelopes''. Now it is customary to call them
``Arens-Michael envelopes'', following the terminology suggested by A.~Ya.~Helemskii \cite{X2}.
As was observed in \cite{Pir_stbflat}, the Arens-Michael envelope of a finitely generated algebra
is a nuclear Fr\'echet algebra.

\begin{example}[\cite{T2}]
\label{ex:AM_poly}
If $A=\CC[x_1,\ldots ,x_n]$ is the polynomial algebra, then $\wh{A}=\cO(\CC^n)$, the
algebra of entire functions on $\CC^n$.
\end{example}

\begin{example}[\cite{Pir_qfree}]
\label{ex:AM_O}
If $(X,\cO^\reg_X)$ is an affine scheme of finite type over $\CC$, and if
$A=\cO^\reg(X)$, then $\wh{A}=\cO(X_h)$, where
$(X_h,\cO_{X_h})$ is the complex space associated to $(X,\cO^\reg_X)$
(cf. \cite[Appendix~B]{Hart_AG}).
\end{example}

Using these examples as a motivation, we defined \cite{Pir_qfree} the algebra $\cO_q(\CC^n)$
of {\em holomorphic
functions on the quantum affine $n$-space} to be the Arens-Michael envelope
of $\cO_q^\reg(\CC^n)$. This algebra can also be described explicitly as follows.
Define a function $w_q\colon \Z_+^n\to\R_+$ by
\begin{equation*}
w_q(k)=
\begin{cases}
1 & \text{if } |q|\ge 1,\\
|q|^{\sum_{i<j} k_i k_j} & \text{if } |q|<1.
\end{cases}
\end{equation*}
As was shown in \cite{Pir_qfree}, we have
\begin{equation}
\label{q_aff_expl}
\cO_q(\CC^n)=
\Bigl\{
a=\sum_{k\in\Z_+^n} c_k x^k :
\| a\|_\rho=\sum_{k\in\Z_+^n} |c_k| w_q(k) \rho^{|k|}<\infty
\;\forall\rho>0
\Bigr\},
\end{equation}
where $|k|=k_1+\cdots +k_n$ for $k=(k_1,\ldots ,k_n)\in\Z_+^n$.
The topology on $\cO_q(\CC^n)$ is given by the norms $\|\cdot\|_\rho \; (\rho>0)$.
Moreover, each norm $\|\cdot\|_\rho$ is submultiplicative.

\section{Quantum polydisk and quantum ball}
The explicit construction \eqref{q_aff_expl} of $\cO_q(\CC^n)$
leads naturally to the following definition.

\begin{definition}[{\cite{Pir_Shlmn,Pir_HFG}}]
Let $q\in\CC^\times$, and let $r>0$.
We define the {\em algebra of holomorphic functions on the
quantum $n$-polydisk of radius $r\in (0,+\infty]$} by
\begin{equation}
\label{q_poly}
\cO_q(\DD^n_r)=
\Bigl\{
a=\sum_{k\in\Z_+^n} c_k x^k :
\| a\|_\rho=\sum_{k\in\Z_+^n} |c_k| w_q(k) \rho^{|k|}<\infty
\;\forall \rho\in (0,r)
\Bigr\}.
\end{equation}
The multiplication on $\cO_q(\DD^n_r)$ is uniquely determined by
$x_i x_j=qx_j x_i\; (i<j)$.
\end{definition}

It follows from the above discussion that $\cO_q(\DD^n_r)$ is a Fr\'echet-Arens-Michael
algebra with respect to the topology determined by the submultiplicative norms
$\|\cdot\|_\rho\; (\rho\in (0,r))$.
It is a simple exercise to show that, if $q=1$, then $\cO_q(\DD^n_r)$ is topologically isomorphic to the algebra
$\cO(\DD^n_r)$ of holomorphic functions on the polydisk
\[
\DD_r^n=\bigl\{ z=(z_1,\ldots ,z_n)\in\CC^n : \max_{1\le i\le n} |z_i|<r\bigr\}.
\]
If $r=\infty$, then we clearly have $\cO_q(\DD^n_r)=\cO_q(\CC^n)$.

The definition of the algebra of holomorphic functions on the quantum ball
is less straightforward. It is based on a theorem by L.~A.~ Aizenberg
and B.~S.~Mityagin \cite{Aiz_Mit}. Recall that
a domain $D\subset\CC^n$ is a \emph{complete Reinhardt domain} if for each
$z=(z_1,\ldots ,z_n)\in D$ and each $(\lambda_1,\ldots ,\lambda_n)\in\CC^n$ satisfying
$|\lambda_i|\le 1\; (i=1,\ldots ,n)$
we have $(\lambda_1 z_1,\ldots ,\lambda_n z_n)\in D$.
Clearly, the polydisk $\DD^n_r$ and the ball
\[
\BB_r^n=\Bigl\{ z=(z_1,\ldots ,z_n)\in\CC^n : \sum_{i=1}^n |z_i|^2<r^2\Bigr\}
\]
are complete Reinhardt domains.

Given a complete bounded Reinhardt domain $D\subset\CC^n$, let
\[
b_k(D)=\sup_{z\in D} |z^k| \qquad (k\in\Z_+^n).
\]
Aizenberg and Mityagin proved that there exists a topological isomorphism
\[
\cO(D)\cong
\Bigl\{
f=\sum_{k\in\Z_+^n} c_k z^k :
\| f\|_s=\sum_{k\in\Z_+^n} |c_k| b_k(D) s^{|k|}<\infty
\;\forall s\in (0,1)\Bigr\}.
\]
Explicitly, the above isomorphism takes each function $f\in\cO(D)$ to its Taylor expansion at $0$.

We clearly have $b_k(\DD^n_r)=r^{|k|}$. An explicit calculation
involving Lagrange's multipliers shows that
\[
b_k(\BB^n_r)=\left(\frac{k^k}{|k|^{|k|}}\right)^{1/2}\!\! r^{|k|}.
\]
Therefore
\begin{equation}
\label{AM_ball_1}
\cO(\BB^n_r)
\cong
\biggl\{
f=\sum_{k\in\Z_+^n} c_k z^k :
\| f\|_\rho=\sum_{k\in\Z_+^n} |c_k| \left(\frac{k^k}{|k|^{|k|}}\right)^{1/2}\!\! \rho^{|k|}<\infty
\;\forall \rho\in (0,r)
\biggr\}.
\end{equation}

Now we have to quantize the above algebra.
To this end, it will be convenient to represent $\cO(\BB^n_r)$ in a slightly different way.

\begin{prop}
\label{prop:ball_powrep}
There exists a topological isomorphism
\begin{equation}
\label{AM_ball_2}
\cO(\BB^n_r)
\cong
\biggl\{
f=\sum_{k\in\Z_+^n} c_k z^k :
\| f\|_\rho=\sum_{k\in\Z_+^n} |c_k| \left(\frac{k!}{|k|!}\right)^{1/2}\!\! \rho^{|k|}<\infty
\;\forall \rho\in (0,r)
\biggr\}.
\end{equation}
\end{prop}

The proof is a simple application of Stirling's formula.

The power series representation \eqref{AM_ball_2} of $\cO(\BB^n_r)$
is more convenient to us than \eqref{AM_ball_1} because there is a standard
way to quantize the factorial. For each $k\in\N$, let
\[
[k]_q=1+q+\cdots +q^{k-1};\quad [k]_q!=[1]_q [2]_q \cdots [k]_q.
\]
It is also convenient to let $[0]_q!=1$. Finally, given $k=(k_1,\ldots ,k_n)\in\Z_+^n$,
we let $[k]_q!=[k_1]_q!\cdots [k_n]_q!$.

\begin{definition}
The \emph{space of holomorphic functions on the quantum $n$-ball}
of radius $r\in (0,+\infty]$ is
\begin{equation}
\label{qball}
\cO_q(\BB^n_r)
=\biggl\{
a=\sum_{k\in\Z_+^n} c_k x^k :
\| a\|_\rho=\sum_{k\in\Z_+^n} |c_k| \left( \frac{[k]_{|q|^{-2}}!}{\bigl[ |k|\bigr]_{|q|^{-2}}!}\right)^{\!\! 1/2}
\!\!\!\! \rho^{|k|}<\infty\; \forall \rho\in (0,r) \biggr\}.
\end{equation}
\end{definition}

It is immediate from the definition that $\cO_q(\BB^n_r)$ is a Fr\'echet space
with respect to the topology determined by the norms $\|\cdot\|_\rho\; (\rho\in (0,r))$.
However, it is not obvious at all whether it is an algebra. In fact it is:

\begin{prop}
The Fr\'echet space $\cO_q(\BB^n_r)$ is an Arens-Michael algebra with respect to
the multiplication uniquely determined by
$x_i x_j=q x_j x_i\quad (i<j)$. Moreover, each norm $\|\cdot\|_\rho$ given by \eqref{qball}
is submultiplicative.
\end{prop}
The proof is based on a $q$-analog of the Chu-Vandermonde formula
(see, e.g., \cite[2.1.2, Proposition 3]{Klim_Schm}).

If $q=1$, then Proposition~\ref{prop:ball_powrep} implies that
$\cO_q(\BB^n_r)\cong \cO(\BB^n_r)$. It can also be shown that, if $r=\infty$,
then $\cO_q(\BB^n_r)=\cO_q(\CC^n)$ (although this is not that obvious as in the case of
the quantum polydisk algebra).

\begin{prop}
\label{prop:q-qinv}
For each $q\in\CC^\times$ and each $r\in (0,+\infty]$ there exist Fr\'echet algebra isomorphisms
\begin{gather*}
\cO_q(\DD^n_r)\cong\cO_{q^{-1}}(\DD^n_r),\quad x_i\mapsto x_{n-i};\\
\cO_q(\BB^n_r)\cong\cO_{q^{-1}}(\BB^n_r),\quad x_i\mapsto x_{n-i}.
\end{gather*}
\end{prop}

The idea of the proof of Proposition \ref{prop:q-qinv} is as follows. Clearly, there is
an algebra isomorphism $\tau\colon\cO_q^\reg(\CC^n)\to\cO_{q^{-1}}^\reg(\CC^n)$
taking each $x_i$ to $x_{n-i}$. An explicit calculation shows that $\tau$ is isometric with respect
to each norm $\|\cdot\|_\rho$, both on $\cO_q(\DD^n_r)$ and $\cO_q(\BB^n_r)$.
The rest is clear.

\section{Quantum ball \`a la Vaksman}
\label{sec:Vaksman}
In the special case where $0<q<1$, $\cO_q(\BB^n_r)$ is closely related to
L.~L.~Vaksman's $q$-analog of $A(\bar\BB^n)$, the algebra of functions
holomorphic on the open unit ball $\BB^n=\BB^n_1$
and continuous on the closed ball $\bar\BB^n$~\cite{Vaks_max}.
Let us recall how Vaksman's algebra is defined.
Assume that $0<q<1$, and let $\Pol_q(\CC^n)$ denote
the $*$-algebra generated (as a $*$-algebra) by $n$ elements $x_1,\ldots ,x_n$
subject to the relations
\begin{equation}
\label{tw_CCR}
\begin{aligned}
x_i x_j&=q x_j x_i \quad (i<j);\\
x_i^* x_j&= q x_j x_i^* \quad (i\ne j);\\
x_i^* x_i&=q^2 x_i x_i^*+(1-q^2)\Bigl(1-\sum_{k>i} x_k x_k^*\Bigr).
\end{aligned}
\end{equation}
Clearly, for $q=1$ we have $\Pol_q(\CC^n)=\Pol(\CC^n)$, where $\Pol(\CC^n)$
is the algebra of polynomials in the complex coordinates $z_1,\ldots ,z_n$ on $\CC^n$
and their complex conjugates $\bar z_1,\ldots ,\bar z_n$.
Observe that the subalgebra of $\Pol_q(\CC^n)$ generated (as an algebra) by $x_1,\ldots ,x_n$
is exactly $\cO_q^\reg(\CC^n)$.
The algebra $\Pol_q(\CC^n)$ was introduced by W.~Pusz and S.~L.~Woronowicz \cite{PW},
although they used different $*$-generators $a_1,\ldots ,a_n$ given by
$a_i=(1-q^2)^{-1/2} x_i^*$. Relations \eqref{tw_CCR} divided by $1-q^2$ and
written in terms of the $a_i$'s
are called the ``twisted canonical commutation relations'', and the algebra $A_q=\Pol_q(\CC^n)$
defined in terms of the $a_i$'s is sometimes called the ``quantum Weyl algebra''
(see, e.g., \cite{WZ,Alev,Jordan,Klim_Schm}).
Note that, while $\Pol_q(\CC^n)$ becomes $\Pol(\CC^n)$ for $q=1$, $A_q$ becomes
the Weyl algebra. The idea to use the generators $x_i$
instead of the $a_i$'s and to consider $\Pol_q(\CC^n)$ as a $q$-analog of $\Pol(\CC^n)$
is probably due to Vaksman \cite{Vaks_splet}; the one-dimensional case
was considered in \cite{Klim_Lesn}.
The algebra $\Pol_q(\CC^n)$ serves as a basic example in the general theory
of quantum bounded symmetric domains developed by Vaksman and his
collaborators (see \cite{Vaks_sborn,Vaks_book} and references therein).

Now let $H$ be a Hilbert space with an orthonormal basis
$\{ e_k : k\in \Z_+^n\}$.
The algebra of bounded linear operators on $H$ will be denoted by $\cB(H)$.
Following \cite{PW}, for each
$k=(k_1,\ldots ,k_n)\in\Z_+^n$ we will write $|k_1,\ldots ,k_n\ra$
for $e_k$. As was proved by Pusz and Woronowicz \cite{PW}, there exists
a faithful irreducible $*$-representation $\pi\colon\Pol_q(\CC^n)\to\cB(H)$ uniquely
determined by
\begin{gather*}
\pi(x_j)e_k=\sqrt{1-q^2} \sqrt{[k_j+1]_{q^2}}\, q^{\sum_{i>j}k_i}
|k_1,\ldots ,k_j+1,\ldots ,k_n\ra\\
(j=1,\ldots ,n,\; k=(k_1,\ldots ,k_n)\in\Z_+^n). \notag
\end{gather*}
The completion of $\Pol_q(\CC^n)$ with respect to the operator norm $\| a\|_\op=\| \pi(a)\|$
is denoted by $C_q(\bar\BB^n)$ and is called the {\em algebra of continuous functions
on the closed quantum ball} \cite{Vaks_max}; see also \cite{PW,Prosk_Sam}.
The closure of $\cO_q^\reg(\CC^n)$ in $C_q(\bar\BB^n)$ is denoted by
$A_q(\bar\BB^n)$ \cite{Vaks_max}; this is a natural $q$-analog of the algebra
$A(\bar\BB^n)$, which consists of those continuous functions on $\bar\BB^n$ that are holomorphic
on $\BB^n$.

For each $\rho>0$, let $\gamma_\rho$ be the automorphism
of $\cO_q^\reg(\CC^n)$ uniquely
determined by $\gamma_\rho(x_i)=\rho x_i\; (i=1,\ldots ,n)$.
Define a norm $\|\cdot\|_{\rho}^\infty$ on $\cO_q^\reg(\CC^n)$ by
\begin{equation*}
\| a\|_{\rho}^\infty=\|\gamma_\rho(a)\|_\op \qquad (a\in\cO_q^\reg(\CC^n)).
\end{equation*}
According to Vaksman's point of view, $\|\cdot\|_\op$ is a natural $q$-analog of
the supremum norm over $\bar\BB^n$. Therefore our $\|\cdot\|_\rho^\infty$
is a $q$-analog of the supremum norm over $\bar\BB^n_\rho$.
It is well known that $\cO^\reg(\CC^n)=\CC[z_1,\ldots ,z_n]$ is dense in $\cO(\DD^n_r)$;
in other words, the completion of $\cO^\reg(\CC^n)$ with respect to the family
$\{\|\cdot\|_\rho^\infty : 0<\rho<r\}$ of norms is topologically isomorphic to $\cO(\BB^n_r)$.
This result has the following $q$-analog.

\begin{theorem}
\label{thm:Vaks_ball}
For each $q\in (0,1)$ and each $r\in (0,+\infty]$,
the completion of $\cO_q^\reg(\CC^n)$ with respect to the family
$\{\|\cdot\|_\rho^\infty : 0<\rho<r\}$ of norms is topologically isomorphic to $\cO_q(\BB^n_r)$.
\end{theorem}

Thus we see that our definition of $\cO_q(\BB^n_r)$ is consistent with the definition
given in \cite{Pir_Shlmn}.

\section{Free polydisk and free ball}
\label{sec:free}
Let $F_n=\CC\la\zeta_1,\ldots ,\zeta_n\ra$ denote the free algebra on $n$ generators
$\zeta_1,\ldots ,\zeta_n$. Clearly, $\cO_q^\reg(\CC^n)$ is nothing but the quotient of $F_n$
modulo the two-sided ideal generated by the elements $\zeta_i\zeta_j-q\zeta_j\zeta_i\; (i<j)$.
Our next goal is to represent the algebras $\cO_q(\DD^n_r)$ and $\cO_q(\BB^n_r)$
in a similar way, i.e., as quotients of certain ``analytic analogs'' of $F_n$. This will enable us to
interpret $\cO_q(\DD^n_r)$ and $\cO_q(\BB^n_r)$ as ``deformations''
(in a suitable sense) of $\cO(\DD^n_r)$ and $\cO(\BB^n_r)$, respectively.

Let us introduce some notation.
For each $d\in \Z_+$, let $W_{n,d}=\{ 1,\ldots ,n\}^d$,
and let $W_n=\bigsqcup_{d\in\Z_+} W_{n,d}$. Thus a typical element of $W_n$
is a $d$-tuple $\alpha=(\alpha_1,\ldots ,\alpha_d)$ of arbitrary length
$d\in\Z_+$, where $\alpha_j\in\{ 1,\ldots ,n\}$ for all $j$.
Given $\alpha\in W_{n,d}\subset W_n$, we let $|\alpha|=d$.
The only element of $W_{n,0}$ will be denoted by $*$.
For each $\alpha=(\alpha_1,\ldots ,\alpha_d)\in W_n$ with $d>0$, let
$\zeta_\alpha=\zeta_{\alpha_1}\cdots\zeta_{\alpha_d}\in F_n$.
It is also convenient to set $\zeta_*=1\in F_n$.
The family $\{\zeta_\alpha : \alpha\in W_n\}$ of all words in $\zeta_1,\ldots ,\zeta_n$
is the standard vector space basis of $F_n$.

Recall from \cite{Pir_HFG} (see also \cite{Pir_Shlmn}) that each family $(A_i)_{i\in I}$
of Arens-Michael algebras has a coproduct $\Pfree_{i\in I} A_i$ in the category
$\ArM$ of Arens-Michael algebras. The algebra $\Pfree_{i\in I} A_i$ is called
the {\em Arens-Michael free product} of the $A_i$'s. By definition \cite{Pir_Shlmn,Pir_HFG},
the {\em algebra of holomorphic functions on the free $n$-polydisk of radius $r\in (0,+\infty]$} is
\[
\cF(\DD^n_r)=\cO(\DD_r)\Pfree\cdots\Pfree\cO(\DD_r).
\]
The algebra $\cF(\DD_r^n)$ can also be described more explicitly as follows.
Let $\fF_n$ denote the algebra of all free formal series $a=\sum_{\alpha\in W_n} c_\alpha\zeta_\alpha$
(where $c_\alpha\in\CC$) with the obvious multiplication.
In other words, $\fF_n=\varprojlim_d F_n/I^d$, where $I$ is the ideal of
$F_n$ generated by $\zeta_1,\ldots ,\zeta_n$.
Given $d\ge 2$ and $\alpha=(\alpha_1,\ldots ,\alpha_d)\in W_{n,d}$, let $s(\alpha)$ denote the
cardinality of the set
\[
\bigl\{ i \in \{ 1,\ldots ,d-1\} : \alpha_i\ne\alpha_{i+1} \bigr\}.
\]
If $d\in\{ 0,1\}$, we let $s(\alpha)=d-1$.
By \cite[Proposition 7.8]{Pir_HFG}, we have
\begin{equation}
\label{free_poly_expl}
\cF(\DD_r^n)=\Bigl\{ a=\sum_{\alpha\in W_n} c_\alpha\zeta_\alpha \in\fF_n :\\
\| a\|_{\rho,\tau}=\sum_{\alpha\in W_n} |c_\alpha|\rho^{|\alpha|} \tau^{s(\alpha)+1}<\infty
\;\forall \rho\in (0,r),\; \forall \tau\ge 1\Bigr\}.
\end{equation}
The topology on $\cF(\DD_r^n)$ is given by the norms
$\|\cdot\|_{\rho,\tau}\; (\rho\in (0,r),\; \tau\ge 1)$,
and the multiplication is given by concatenation.
Moreover, each norm $\|\cdot\|_{\rho,\tau}$ is submultiplicative.

Another natural candidate for the algebra of holomorphic functions
on the free polydisk was introduced by J.~L.~Taylor \cite{T2,T3}.
By definition,
\begin{equation}
\label{F_poly_T}
\cF^T(\DD_r^n)=\Bigl\{ a=\sum_{\alpha\in W_n} c_\alpha\zeta_\alpha \in\fF_n :
\| a\|_\rho=\sum_{\alpha\in W_n} |c_\alpha|\rho^{|\alpha|}<\infty
\;\forall \rho\in (0,r)\Bigr\}.
\end{equation}
It is easy to see that $\cF(\DD^n_r)\subset\cF^T(\DD^n_r)$, and that the
inclusion of $\cF(\DD^n_r)$ into $\cF^T(\DD^n_r)$ is continuous.
On the other hand, $\cF(\DD^n_r)\ne\cF^T(\DD^n_r)$ unless $n=1$ or $r=\infty$.
Note also that, if $r=\infty$, then both $\cF(\DD^n_r)$ and $\cF^T(\DD^n_r)$ are
topologically isomorphic to $\cF(\CC^n)$, the Arens-Michael envelope of $F_n$.

\begin{theorem}
\label{thm:q_poly_quot_free_poly}
Let $q\in\CC^\times$, $n\in\N$, and $r\in (0,+\infty]$.
\begin{compactenum}
\item[{\upshape (i)}]
The algebra $\cO_q(\DD^n_r)$ is topologically isomorphic to the quotient
of $\cF(\DD^n_r)$ modulo the closed two-sided ideal generated by the elements
$\zeta_i\zeta_j-q\zeta_j\zeta_i\; (i<j)$. Moreover, for each $\rho\in (0,r)$ and each $\tau\ge 1$
the norm $\|\cdot\|_\rho$ on $\cO_q(\DD^n_r)$ given by \eqref{q_poly} is equal to
the quotient of the norm $\|\cdot\|_{\rho,\tau}$ on $\cF(\DD^n_r)$ given by \eqref{free_poly_expl}.
\item[{\upshape (ii)}]
The algebra $\cO_q(\DD^n_r)$ is topologically isomorphic to the quotient
of $\cF^T(\DD^n_r)$ modulo the closed two-sided ideal generated by the elements
$\zeta_i\zeta_j-q\zeta_j\zeta_i\; (i<j)$. Moreover, for each $\rho\in (0,r)$
the norm $\|\cdot\|_\rho$ on $\cO_q(\DD^n_r)$ given by \eqref{q_poly} is equal to
the quotient of the norm $\|\cdot\|_\rho$ on $\cF^T(\DD^n_r)$ given by \eqref{F_poly_T}.
\end{compactenum}
\end{theorem}

Part (i) of Theorem \ref{thm:q_poly_quot_free_poly}, except for the
equality of the norms, was proved in \cite[Theorem 7.13]{Pir_HFG}
in the more general multiparameter case. Part (ii) is new.

To formulate a similar result for $\cO_q(\BB^n_r)$, we need G. Popescu's algebra of
``holomorphic functions on the free ball'' \cite{Pop_holball}.
Let $H$ be a Hilbert space, and
let $T=(T_1,\ldots ,T_n)$ be an $n$-tuple of bounded linear operators on $H$.
Following \cite{Pop_holball}, we identify $T$ with the ``row'' operator
acting from the Hilbert direct sum $H^n=H\oplus\cdots\oplus H$ to $H$.
Thus we have $\| T\|=\|\sum_{i=1}^n T_i T_i^*\|^{1/2}$.
For each $a=\sum_\alpha c_\alpha\zeta_\alpha\in\fF_n$,
the {\em radius of convergence} $R(a)\in [0,+\infty]$
is given by the Cauchy-Hadamard-type formula
\[
\frac{1}{R(a)}=\limsup_{d\to\infty}%
\Bigl(\sum_{|\alpha|=d} |c_\alpha|^2\Bigr)^{\frac{1}{2d}}.
\]
By \cite[Theorem 1.1]{Pop_holball}, for each $T\in\cB(H)^n$ such that
$\| T\|<R(a)$, the series
\begin{equation}
\label{Pop_series}
\sum_{d=0}^\infty \Bigl(\sum_{|\alpha|=d} c_\alpha T_\alpha \Bigr)
\end{equation}
converges in $\cB(H)$ and, moreover, $\sum_d \| \sum_{|\alpha|=d} c_\alpha T_\alpha\|<\infty$.
On the other hand, if $H$ is infinite-dimensional,
then for each $R'>R(a)$ there exists $T\in\cB(H)^n$ with $\| T\|=R'$
such that the series~\eqref{Pop_series} diverges. This ``free operator analog'' of the
classical Hadamard lemma explains why the radius of convergence
is so called.

By \cite[Theorem 1.4]{Pop_holball}, the collection of all $a\in\fF_n$ such that
$R(a)\ge r$ is a subalgebra of $\fF_n$. We denote this algebra by $\cF(\BB^n_r)$
(Popescu uses the notation $\cHol(\cB(\ccH)_r^n)$),
and we call it the {\em algebra of holomorphic functions on the free $n$-ball of radius $r$} \cite{Pir_Shlmn}.
For each $a\in \cF(\BB^n_r)$, each Hilbert space $H$,
and each $T\in\cB(H)^n$ with $\| T\|<r$, the sum of the series~\eqref{Pop_series}
is denoted by $a(T)$. The map
\[
\cF(\BB^n_r) \to \cB(H),\qquad a\mapsto a(T),
\]
is an algebra homomorphism.

Fix an infinite-dimensional Hilbert space $\ccH$, and,
for each $\rho\in (0,r)$, define a seminorm $\|\cdot\|_\rho^P$ on $\cF(\BB^n_r)$
by
\[
\| a\|_\rho^P=\sup\{ \| a(T)\| : T\in\cB(\ccH)^n,\; \| T\|\le\rho\}.
\]
As was observed in \cite{Pop_holball}, $\|\cdot\|_\rho^P$ is in fact a norm on $\cF(\BB^n_r)$.
This norm can be viewed as a ``free operator analog'' of the supremum norm over $\bar\BB^n_\rho$.
By \cite[Theorem 5.6]{Pop_holball}, $\cF(\BB^n_r)$ is a Fr\'echet algebra with respect
to the topology determined by the family $\{ \|\cdot\|_\rho^P : \rho\in (0,r)\}$ of norms.

For our purposes, a slightly different definition of $\cF(\BB^n_r)$ is needed.
Consider the projection $p\colon W_n\to\Z_+^n$ given by
\[
p(\alpha)=(p_1(\alpha),\ldots ,p_n(\alpha)),\quad
p_j(\alpha)=|\alpha^{-1}(j)|.
\]

\begin{prop}
\label{prop:freeball}
There exists a topological isomorphism
\begin{equation}
\label{freeball}
\cF(\BB_r^n)\cong\biggl\{ a=\sum_{\alpha\in W_n} c_\alpha\zeta_\alpha\in\fF_n :
\| a\|_\rho=\sum_{k\in\Z_+^n}\Bigl(\sum_{\alpha\in p^{-1}(k)}%
|c_\alpha|^2\Bigr)^{1/2} \rho^{|k|}<\infty
\;\forall \rho\in (0,r)\biggr\}.
\end{equation}
Moreover, each norm $\|\cdot\|_\rho$ given by \eqref{freeball} is submultiplicative.
\end{prop}

Using Proposition \ref{prop:freeball}, we obtain the following theorem, which extends
our earlier result from \cite{Pir_Shlmn}.

\begin{theorem}
\label{thm:q_ball_quot_free_ball}
For each $q\in\CC^\times$, the algebra $\cO_q(\BB^n_r)$ is topologically isomorphic to the quotient
of $\cF(\BB^n_r)$ modulo the closed two-sided ideal generated by the elements
$\zeta_i\zeta_j-q\zeta_j\zeta_i\; (i<j)$. Moreover, for each $\rho\in (0,r)$
the norm $\|\cdot\|_\rho$ on $\cO_q(\BB^n_r)$ given by \eqref{qball} is equal to
the quotient of the norm $\|\cdot\|_\rho$ on $\cF(\BB^n_r)$ given by \eqref{freeball}.
\end{theorem}

\section{Fr\'echet algebra bundles}

Now we can explain in which sense the algebras $\cO_q(\DD^n_r)$ and $\cO_q(\BB^n_r)$
are ``deformations'' of $\cO(\DD^n_r)$ and $\cO(\BB^n_r)$, respectively.
Let us recall some definitions from \cite{Gierz} (in a slightly modified form).
Suppose that $X$ is a locally compact, Hausdorff topological space.
By a {\em family of vector spaces} over $X$ we mean a pair $(E,p)$, where
$E$ is a set and $p\colon E\to X$ is a surjective map, together with a vector space structure
on each fiber $E_x=p^{-1}(x)\; (x\in X)$. As usual, we let
\[
E\times_X E=\{ (u,v)\in E\times E : p(u)=p(v)\}.
\]
By a {\em prebundle of topological vector spaces}
over $X$ we mean a family $(E,p)$ of vector spaces over $X$ together with a topology on $E$
such that $p$ is continuous and open, the zero section $0\colon X\to E$ is continuous,
and the operations
\begin{equation*}
\begin{split}
E\times_X E&\to E, \quad (u,v)\mapsto u+v,\\
\CC\times E&\to E,\quad (\lambda, u)\mapsto \lambda u,
\end{split}
\end{equation*}
are also continuous.

Let $(E,p)$ be a family of vector spaces over $X$.
By definition, a function $\|\cdot\|\colon E\to [0,+\infty)$ is a {\em seminorm}
if the restriction of $\|\cdot\|$ to each fiber is a seminorm in the usual sense.
A family $\cN=\{ \|\cdot\|_\lambda : \lambda\in\Lambda\}$ of seminorms on $E$ is said to be
{\em directed} if for each $\lambda,\mu\in\Lambda$ there exist $C>0$ and $\nu\in\Lambda$
such that $\|\cdot\|_\lambda\le C\|\cdot\|_\nu$ and $\|\cdot\|_\mu\le C\|\cdot\|_\nu$.
If $\cN=\{ \|\cdot\|_\lambda : \lambda\in\Lambda\}$ and
$\cN'=\{ \|\cdot\|_\mu : \mu\in\Lambda'\}$ are two directed families of seminorms on $E$,
then we say that $\cN$ is {\em dominated} by $\cN'$ (and write $\cN\prec\cN'$)
if for each $\lambda\in\Lambda$ there exist $C>0$ and $\mu\in\Lambda'$
such that $\|\cdot\|_\lambda\le C\|\cdot\|_\mu$.
If $\cN\prec\cN'$ and $\cN'\prec\cN$, then we say that $\cN$ and $\cN'$ are
{\em equivalent} and write $\cN\sim\cN'$.

\begin{remark}
\label{rem:uniform_str}
A directed family $\cN=\{ \|\cdot\|_\lambda : \lambda\in\Lambda\}$ of seminorms on $E$
determines a uniform structure $\cU(\cN)$ on $E$ whose basis consists of all sets of the form
\[
\bigl\{ (u,v)\in E\times_X E : \| u-v\|_\lambda<\eps\bigr\}
\qquad (\lambda\in\Lambda,\; \eps>0).
\]
It is easy to see that $\cN\prec\cN'$ if and only if $\cU(\cN)\subset\cU(\cN')$, and
consequently $\cN\sim\cN'$ if and only if $\cU(\cN)=\cU(\cN')$.
\end{remark}

The following definition is essentially a locally convex version of the notion of a Banach bundle
in the sense of J.~M.~G.~Fell (see, e.g., \cite{FD1}).

\begin{definition}
Let $(E,p)$ be a prebundle of topological vector spaces over $X$, and let
$\cN=\{ \|\cdot\|_\lambda : \lambda\in\Lambda\}$ be a directed family of seminorms on $E$.
We say that $\cN$ is {\em admissible} if for each $x\in X$ the sets
\[
\{ u\in E : p(u)\in U,\; \| u\|_\lambda<\eps\}\qquad
(\lambda\in\Lambda,\; \eps>0,\; U\subseteq X\text{ is an open neighborhood of }x)
\]
form a base of open neighborhoods of $0\in E_x$.
By a {\em locally convex uniform structure} on $(E,p)$ we mean the
equivalence class of an admissible directed family of seminorms on $E$.
By a {\em locally convex bundle} over $X$ we mean a
prebundle of topological vector spaces over $X$ together with a locally convex uniform structure.
A {\em locally convex algebra bundle} over $X$ is a locally
convex bundle $(E,p)$ over $X$ together
with an algebra structure on each fiber $E_x$ such that the map
\begin{equation*}
m_E\colon E\times_X E\to E, \quad (u,v)\mapsto uv,
\end{equation*}
is continuous. A {\em Fr\'echet algebra bundle} over $X$
is a locally convex algebra bundle $(E,p)$ over $X$
such that each fiber $E_x$ is a Fr\'echet algebra.
\end{definition}

The following result is an improvement of \cite[Theorem 6.1]{Pir_Shlmn}.

\begin{theorem}
\label{thm:bnd}
Let $n\in\N$, and let $r\in (0,+\infty]$.
\begin{compactenum}
\item[{\upshape (i)}] There exists a Fr\'echet algebra bundle
$(D,p)$ over $\CC^\times$ such that for each $q\in\CC^\times$ we have
$D_q\cong\cO_q(\DD_r^n)$.
\item[{\upshape (ii)}] There exists a Fr\'echet algebra bundle
$(B,p)$ over $\CC^\times$ such that for each $q\in\CC^\times$ we have
$B_q\cong\cO_q(\BB^n_r)$.
\end{compactenum}
\end{theorem}

In fact, we constructed the bundles $(D,p)$ and $(B,p)$ already in \cite[Theorem 6.1]{Pir_Shlmn},
but we did not know how the fibers of $(B,p)$ look like unless $q\in (0,1]$.
The construction is as follows. Let $z$ denote the complex coordinate on $\CC^\times$,
and let $\cO_\defo(\DD^n_r)$ denote the quotient of $\cO(\CC^\times,\cF(\DD^n_r))$
modulo the closed two-sided ideal generated by the elements $x_i x_j-zx_j x_i\; (i<j)$.
The algebra $\cO_\defo(\DD^n_r)$ is a Fr\'echet $\cO(\CC^\times)$-algebra in a canonical way
(i.e., we have a continuous homomorphism from $\cO(\CC^\times)$ to the center of
$\cO_\defo(\DD^n_r)$). Given a Fr\'echet $\cO(\CC^\times)$-algebra $R$, we define the
{\em fiber} of $R$ over $q\in\CC^\times$ to be $R_q=R/\ol{M_q R}$, where
$M_q=\{ f\in\cO(\CC^\times) : f(q)=0\}$. Let now $E=\bigsqcup_{q\in\CC^\times} R_q$,
and let $p\colon E\to\CC^\times$ take each $R_q$ to $q$.
There is a canonical way to topologize $E$ in such a way that $(E,p)$ becomes
a Fr\'echet algebra bundle over $\CC^\times$. Applying this construction to $R=\cO_\defo(\DD^n_r)$,
we obtain a Fr\'echet algebra bundle $(D,p)$ over $\CC^\times$ whose fibers are equal to those
of $\cO_\defo(\DD^n_r)$. On the other hand, Theorem~\ref{thm:q_poly_quot_free_poly} (i)
implies that the fiber of $\cO_\defo(\DD^n_r)$ over $q\in\CC^\times$ is isomorphic to $\cO_q(\DD^n_r)$.
We can also use $\cF^T(\DD^n_r)$ instead of $\cF(\DD^n_r)$ in the above construction;
the deformation algebra $\cO_\defo(\DD^n_r)$ and hence the bundle $(D,p)$ will then be the same.
By replacing $\cF(\DD^n_r)$ with $\cF(\BB^n_r)$ is the above construction, and by using
Theorem~\ref{thm:q_ball_quot_free_ball} instead of Theorem~\ref{thm:q_poly_quot_free_poly},
we obtain a Fr\'echet $\cO(\CC^\times)$-algebra $\cO_\defo(\BB^n_r)$ and a Fr\'echet
algebra bundle $(B,p)$ whose fiber over $q\in\CC^\times$ is isomorphic to $\cO_q(\BB^n_r)$.

\begin{remark}
We hope that Theorem \ref{thm:bnd} may be of interest from the viewpoint of ``analytic''
deformation theory of associative algebras. Recall that
most papers and monographs on algebraic deformation theory deal with \emph{formal}
deformations. Roughly, this means that the deformed product of two elements of an algebra $A$
is no longer an element of $A$, but is a formal power series with coefficients in $A$.
To some extent, such an approach to deformation
theory is dictated by convenience considerations. However, formal deformations are not
entirely satisfactory from the point of view of physics. Indeed, only those deformations
are physically interesting that are represented by convergent power series \cite{conv_wick}.
In other words, one should be able to substitute concrete values of Planck's constant
into the deformed product. The first successful attempt to develop a theory
satisfying this requirement was made by M.~Rieffel \cite{Rf_Heis,Rf_dq_oa,Rf_Lie,Rf_mem}.
In his approach, a deformation is a continuous field (or a bundle) of $C^*$-algebras endowed with
an additional structure. See \cite{Rf_quest} for a recent survey.

For a number of reasons, it is natural to expect that there should be an approach
intermediate between the formal and continuous (i.e., $C^*$-algebraic) deformation theories.
Specifically, it is desirable to have a ``differentiable'' or ``complex analytic'' deformation
theory. Perhaps the first attempt to develop such a theory was made
by M.~J.~Pflaum and M.~Schottenloher \cite{Pfl_Schott}. Important special cases of convergent
deformed products were studied by
H.~Omori, Y.~Maeda, N.~Miyazaki, and A.~Yoshioka \cite{OMMY1,OMMY2}, and by S.~Beiser,
H.~R\"omer, and S.~Waldmann \cite{conv_wick}. Quite recently, very interesting preprints by
S.~Beiser, G.~Lechner, and S.~Waldmann \cite{BW,LW} have appeared, in which a rather general
approach to analytic deformation theory has been presented.
On the other hand, so far nothing is known about deformations of algebras of holomorphic functions
on classical domains. We do not try to give a general definition of a Fr\'echet algebra deformation here,
but we believe that such a definition should be similar to that given by Rieffel in the $C^*$-algebra case.
Also, we would like to note that our deformation algebras $\cO_\defo(\DD^n_r)$
and $\cO_\defo(\BB^n_r)$ are not topologically free (moreover, they are not topologically
projective) over $\cO(\CC^\times)$, so they do not fit into the framework suggested in \cite{Pfl_Schott}.
\end{remark}

\section{A $q$-analog of Poincar\'e's theorem}

We now turn to our main question of whether $\cO_q(\DD^n_r)$ and $\cO_q(\BB^n_r)$
are topologically isomorphic. The answer is as follows.

\begin{theorem}
\label{thm:q-Poincare}
Let $n\in\N$ and $r\in (0,+\infty]$.
\begin{compactenum}
\item[{\upshape (i)}]
If $n\ge 2$, $r<\infty$, and $|q|=1$, then $\cO_q(\DD^n_r)$
and $\cO_q(\BB^n_r)$ are not topologically isomorphic.
\item[{\upshape (ii)}]
If $|q|\ne 1$, then
$\cO_q(\DD^n_r)$ and $\cO_q(\BB^n_r)$ are topologically isomorphic
(in fact, they are equal as power series algebras).
\end{compactenum}
\end{theorem}

The proof of part (ii) is elementary in the case where $|q|>1$. The case where $|q|<1$
is reduced to the previous one via Proposition~\ref{prop:q-qinv}.

The proof of part (i) is more involved and is based on an Arens-Michael algebra version of
the joint spectral radius (see, e.g., \cite{Muller_book}).
Let $A$ be an Arens-Michael algebra, and let $\{ \|\cdot\|_\lambda : \lambda\in\Lambda\}$
be a directed defining family of submultiplicative seminorms on $A$.
Given an $n$-tuple $a=(a_1,\ldots ,a_n)\in A^n$, we define the {\em joint $\ell^p$-spectral radius}
$r^A_p(a)$ by
\begin{equation*}
\begin{split}
r^A_p(a)&=\sup_{\lambda\in\Lambda}%
\lim_{d\to\infty}\Bigl(\sum_{\alpha\in W_{n,d}} \| a_\alpha\|_\lambda^p\Bigr)^{1/pd}
\quad\text{for }1\le p<\infty;\\
r^A_\infty(a)&=\sup_{\lambda\in\Lambda}%
\lim_{d\to\infty}\Bigl(\sup_{\alpha\in W_{n,d}} \| a_\alpha\|_\lambda\Bigr)^{1/d}.
\end{split}
\end{equation*}
In the case of Banach algebras, the joint $\ell^\infty$-spectral radius was introduced
by G.-C.~Rota and W.~G.~Strang \cite{Rota_Str}.
The case $p<\infty$ was studied by A.~So\l tysiak \cite{Sol_jsprad}
for commuting $n$-tuples in a Banach algebra. The observation that a similar definition
makes sense in the noncommutative case is probably due to V.~M\"uller \cite[C.32.2]{Muller_book}.

It is easy to show that $r^A_p(a)$ does not depend on the choice of a
directed defining family of submultiplicative seminorms on $A$.
Moreover, if $\varphi\colon A\to B$ is a continuous homomorphism of Arens-Michael algebras,
then for each $a=(a_1,\ldots ,a_n)\in A^n$ we have $r^B_p(\varphi(a))\le r^A_p(a)$,
where $\varphi(a)=(\varphi(a_1),\ldots ,\varphi(a_n))$.
In particular, if $\varphi\colon A\to B$ is a topological isomorphism, then
$r^B_p(\varphi(a))=r^A_p(a)$.
In contrast to the Banach algebra case, it may happen that $r^A_p(a)=+\infty$.
For example, if $A=\cO(\CC)$ and $z\in A$ is the complex coordinate, then an easy
computation shows that $r_p^A(z)=+\infty$ for all $p$.

The following example is crucial for our purposes.
Let $|q|=1$, and let $x=(x_1,\ldots ,x_n)$ denote the system of canonical
generators of $\cO_q(\DD^n_r)$ or $\cO_q(\BB^n_r)$. Then we have
\begin{equation}
\label{jsprad_qpoly_qball}
r_2^{\cO_q(\BB^n_r)}(x)=r\quad\text{and}\quad
r_2^{\cO_q(\DD^n_r)}(x)=r\sqrt{n}.
\end{equation}

Now we can explain the idea of the proof of Theorem~\ref{thm:q-Poincare} (i).
Fix $q\in\CC^\times,\; q\ne 1$, with $|q|=1$.
Let $A=\cO_q(\DD^n_r)$ and $B=\cO_q(\BB^n_r)$, and
assume that $\varphi\colon B\to A$ is a topological isomorphism.
For each $i=1,\ldots ,n$, let $f_i=\varphi(x_i)$. A tedious but elementary algebraic argument
shows that there exists a permutation $\sigma\in S_n$ such that
\[
f_i=\lambda_i x_{\sigma(i)}+\text{terms of higher degree},
\]
where $|\lambda_i|=1$. Hence for each $\alpha=(\alpha_1,\ldots ,\alpha_d)\in W_n$ we have
\[
f_\alpha=\lambda_\alpha x_{\sigma(\alpha)}+\text{terms of higher degree},
\]
where $\sigma(\alpha)=(\sigma(\alpha_1),\ldots ,\sigma(\alpha_d))$ and
$\lambda_\alpha=\lambda_{\alpha_1}\cdots\lambda_{\alpha_d}$.
This implies that for each $\rho\in (0,r)$ we have
$\| f_\alpha\|_\rho \ge \| x_{\sigma(\alpha)}\|_\rho=\| x_\alpha\|_\rho$,
where $\|\cdot\|_\rho$ is the norm on $\cO_q(\DD^n_r)$ given by~\eqref{q_poly}.
Hence $r_2^A(f)\ge r_2^A(x)$. On the other hand, $r_2^A(f)=r_2^B(x)$, because
$\varphi$ is a topological isomorphism. Taking into account \eqref{jsprad_qpoly_qball},
we conclude that $r=r_2^B(x)\ge r_2^A(x)=r\sqrt{n}$, whence $n=1$.

\section{Open problems}
We conclude the paper with a couple of open problems. The first problem was already
mentioned in Section~\ref{sec:Vaksman} and is inspired by the fact that
the algebra $\cO_q(\BB^n_r)$ is defined for all $q\in\CC^\times$.

\begin{problem}
Is it possible to define Vaksman's algebras $C_q(\bar\BB^n)$ and $A_q(\bar\BB^n)$
in the case where $q\notin (0,1]$? If yes, then can Theorem~{\upshape\ref{thm:Vaks_ball}}
be extended to this case?
\end{problem}

The second problem is related to the notion of an HFG algebra introduced in \cite{Pir_Shlmn,Pir_HFG}.
Let $\cF(\CC^n)$ denote the Arens-Michael envelope of the free algebra $F_n$
(see Section~\ref{sec:free}).
A Fr\'echet algebra $A$ is said to be {\em holomorphically finitely generated} (HFG for short)
if $A$ is isomorphic to a quotient of $\cF(\CC^n)$ for some $n$.
There is also an ``internal'' definition given in terms of J.~L.~Taylor's free functional
calculus \cite{T2}. By \cite[Theorem 3.22]{Pir_HFG},
a commutative Fr\'echet-Arens-Michael algebra is holomorphically finitely generated if and only if
it is topologically isomorphic to $\cO(X)$ for some Stein space $(X,\cO_X)$ of finite
embedding dimension. Together with Forster's theorem (see Section~\ref{sec:intro}),
this implies that the category of commutative HFG algebras is anti-equivalent to the
category of Stein spaces of finite embedding dimension. There are many natural examples
of noncommutative HFG algebras; see \cite[Section 7]{Pir_HFG}.
For instance, $\cO_q(\DD^n_r)$ and $\cF(\DD^n_r)$ are HFG algebras.
By Theorem~\ref{thm:q-Poincare} (ii), $\cO_q(\BB^n_r)$ is an HFG algebra
provided that $|q|\ne 1$.

\begin{problem}
Is $\cO_q(\BB^n_r)$ an HFG algebra in the case where $|q|=1,\; q\ne 1$?
\end{problem}

\begin{ackn}
The author thanks A.~Ya.~Helemskii and D. Proskurin for helpful discussions.
\end{ackn}

\end{document}